\input amssym.def
\input amssym
\magnification=1200
\parindent0pt
\hsize=16 true cm
\baselineskip=13  pt plus .2pt
$ $

\def\G{(\Gamma,{\cal G})}
\def\Z{{\Bbb Z}}
\def\D{{\Bbb D}}
\def\A{{\Bbb A}}
\def\S{{\Bbb S}}

\centerline {\bf  On upper and lower bounds for finite group actions on bounded surfaces,}

\centerline {\bf  handlebodies, closed handles and finite graphs}

\bigskip

\centerline {Bruno P. Zimmermann}

\medskip

\centerline {Universit\`a degli Studi di Trieste}

\centerline {Dipartimento di Matematica e Geoscienze}

\centerline {34127 Trieste, Italy}

\bigskip  \bigskip

{\it Abstract.}  This is a survey on upper and lower bounds for finite
group actions on bounded surfaces, 3-dimensional handlebodies and closed handles,
handlebodies in arbitrary dimensions and finite graphs (the common feature of these objects
is that all have free fundamental group).

\bigskip \bigskip

{\bf 1.  Presentation of results}

\medskip

We consider finite group actions of large orders on various low-dimensional manifolds with
free fundamental group, and also on finite graphs. All group actions in the present paper
will faithful and smooth.

\bigskip

{\bf 1.1}  \hskip 1mm  Finite group actions on 3-dimensional handlebodies and bounded
surfaces.

\medskip

In analogy with the classical Hurwitz-bound
$84(g-1)$ for the order of a finite, orienta-tion-preserving group action on a  closed,
orientable surface of genus
$g \ge 2$, an uper bound for the order of a finite group of orientation-preserving
diffeomorphisms of a a 3-dimension handlebody $V_g$ of genus $g \ge 2$ is $12(g-1)$
([Z5],[McMZ1, Theorem 7.2]). More generally, the following holds.

\bigskip

{\bf Theorem 1.} ([MiZ])  {\sl   Let $m_{hb}(g)$ denote the maximum order of a finite,
orientation-preserving group action on a 3-dimensional handlebody of genus $g > 1$.

\medskip

i) There are the upper and lower bounds  
$$4(g+1) \;  \le \;\;   m_{hb}(g) \;\;  \le \; 12(g-1),$$ 
and both $4(g+1)$ and  $12(g-1)$ occur for infinitely many genera $g$.

\medskip

ii) If $g$ is odd then  $m_{hb}(g) = 8(g-1)$ or $\; m_{hb}(g) = 12(g-1)$, and both cases
occur  for infinitely many values of $g$.

\medskip

iii)  The possible values of $m_{hb}(g)$ are of the form ${4n \over n-2}(g-1)$, for an
integer $n \ge 3$, and infinitely many values of $n$ occur resp. do not occur. Moreover if
a value of $n$ occurs then it occurs for infinitely many $g$. }

\bigskip

For finite group actions on bounded surfaces (compact with
nonempty boundary, orientable or not), exactly the same results hold (and with the same
proofs), using the setting in [McMZ2] (see [MiZ, section 3]). Note that, by taking the
product with an interval (twisted if the surface is nonorientable), every finite group
acting on a bounded surface of {\it algebraic genus} $g$ (defined as the rank of the free
fundamental group)  admits also an orientation-preserving action on a handlebody of genus
$g$.

\bigskip

{\bf Theorem 2.}  {\sl  Let $\; m_{bs}(g)$ denote the maximum order of a finite, possibly
orientation-reversing group action on a bounded, orientable or nonorientable surface of
algebraic genus  $g > 1$.  

\medskip

i)  $\;\;  m_{bs}(g) \; \le  \; m_{hb}(g)$

\medskip

ii)  All statements of Theorem 1 remain true for $\; m_{bs}(g)$. 

\medskip

iii) There are values of $g$ such that $\; m_{bs}(g)$ is strictly smaller than 
$\; m_{hb}(g)$.}

\bigskip

Part iii) of Theorem 2 is proved in [CZ] by computational methods;  
the two smallest values of $g$ such that $m_{bs}(g) <  m_{hb}(g)$ are $g = 161$ and $g =
3761$.

\bigskip

{\bf 1.2}  \hskip 1mm  Finite group actions on closed handles.

\medskip

After the classical cases of 3-dimensional handlebodies and bounded surfaces, we consider
actions of finite groups $G$ on closed 3-dimensional analogues of handlebodies, the
connected sums  $H_g = \sharp_g (S^1 \times S^2)$ of $g$ copies of $S^1 \times S^2$
(similar as 
$V_g = \sharp_g^\partial (S^1 \times D^2)$ is the boundary-connected sum of $g$ copies of
$S^1 \times D^2$; so $H_g$ is the double of $V_g$ along its boundary).  We will call $H_g$
a {\it closed handle} or just a {\it handle} of genus $g$

\medskip

Since $H_g$  admits $S^1$-actions (see [R]), it admits finite cyclic group actions of
arbitrarily large order acting trivially on the fundamental group. Let
$G_0$ denote the normal subgroup of all elements of $G$ acting trivially on the
fundamental group (up to inner automorphisms); by [Z1, Proposition 2], 
$G_0$ is cyclic, the quotient $H_g/G_0$ is again a closed handle
of the same genus $g$ and the factor group
$G/G_0$ acts faithfully on the fundamental group of the quotient $H_g/G_0 \cong H_g$. Hence
one is led to consider actions of finite groups $G$ on
$H_g$ which act faithfully on the fundamental group, i.e. induce an
injection into the outer automorphism group ${\rm Out} \, F_g$ of the fundamental group of
$H_g$, the free group $F_g$ of rank $g$.

\bigskip

{\bf Theorem 3.}  ([Z1])  {\sl  Let $m_{ch}(g)$ denote the maximum order of a finite,
orientation-preserving group action on a closed handle $H_g$ of genus $g > 1$ which
induces a faithful action on the fundamental group.

\medskip

i)  For $g \ge 15$, there is the quadratic upper bound  $\; m_{ch}(g) \le  \; 24g(g-1)$.

\medskip

ii) For all $g$, there are the quadratic lower bounds  
$\;  2 g^2  \le  m_{ch}(g) \;$  if $g$ is even, and  
$\; (g+1)^2  \le  \; m_{ch}(g) \;$ if $g$ is odd.}

\bigskip

We don't know the exact value of $\; m_{ch}(g)$ at present but believe  that for large $g$
it coincides with the lower bounds $\; 2 g^2$ resp. $\; (g+1)^2$ of the second part of
Theorem 3; for small values of $g$ there are group actions of larger orders, e.g. 
$m_{ch}(2) = 12$, $\; m_{ch}(3) = 48$ and $\; m_{ch}(4) = 192$.

\bigskip

Next we consider the case of {\it free} actions of finite groups on closed handles $H_g$
which is in  strong analogy with the cases of arbitrary (i.e., not necessarily free) actions
on handlebodies and bounded surfaces  (where free means that every nontrivial element has
empty fixed point set).

\bigskip

{\bf Theorem 4.}  {\sl  Let  $m(g) = \bar m_{ch}(g)$ denote the maximum order of a
free, orientation-preserving finite group action on a closed handle $H_g$ of genus
$g>1$. 

\medskip

i) For all $g > 1$,  
$$2(g+1) \;  \le \;  m(g) \; \le \; 6(g-1)$$ 
and both $2(g+1)$ and  $6(g-1)$ occur for infinitely many genera $g$.

\medskip

ii) If $g$ is odd then  $\; m(g) = 4(g-1)$ or $\; m(g) = 6(g-1)$, and both
cases occur for infinitely many $g$.

\medskip

iii)  The possible values of $m(g)$ are of the form ${2n \over n-2}(g-1)$, for
an integer $n \ge 3$, and infinitely many values of $n$ occur resp. do not occur.}

\bigskip

We note that exactly the same results hold for finite, orientation-preserving group actions
on bounded, orientable surfaces of algebraic genus $g$.

\medskip

The proof of Theorem 4 combines  methods of the handlebody case (Theorem 1) with those for
closed handles (Theorem 3); since it is shorter and less technical, as
an illustration of the methods we give the proof in section 2.

\bigskip

{\bf 1.3}  \hskip 1mm  Finite group actions on handlebodies and closed handles in arbitrary
dimensions.

\medskip

A closed handle $H_g$ is the boundary of a 4-dimensional handlebody, in
particular Theorem 3 holds also for finite group actions on 4-dimensional handlebodies. More
generally,  an orientable handlebody $V^d_g$ of dimension $d$ and genus $g$ is defined 
as a regular neighbourhood of a finite graph, with free fundamental group of
rank $g$, embedded in the sphere $S^d$; alternatively, it is obtained from the closed
disk $D^d$ of dimension $d$ by attaching along its boundary  $g$ copies of a handle $D^{d-1}
\times [0,1]$ in an orientable way, or as the boundary-connected sum 
$\; \sharp_g^\partial (S^1 \times D^{d-1})$ of $g$
copies of $S^1 \times D^{d-1}$. The boundary of $V^d_g$ is a closed manifold
$H^{d-1}_g$ which is the connected sum  $\; \sharp_g (S^1 \times S^{d-2})$  of $g$ copies of
$S^1 \times S^{d-2}$, so these are the higher-dimensional analogues of the closed handles
$H_g = H^3_g$.

\medskip

After Thurston  and Perelman, finite group actions in dimension 3 are geometric; this is no
longer true in higher dimensions, so in order to generalize Theorem 3 one has to consider
some kind of standard actions also in higher dimensions. A natural way to proceed is to
uniformize handlebodies $V^d_g$ by Schottky groups (free groups of M\"obius
transformations of $D^d$ acting by isometries on its interior, the 
Poincar\'e-model of hyperbolic space $\Bbb H^d$); this realizes the interior of a
handlebody $V^d_g$ as a hyperbolic manifold, and we will consider finite groups of
isometries of such hyperbolic (Schottky type) handlebodies.

\medskip

By [Z2], every finite subgroup of the outer automorphism group  ${\rm Out} \, F_g$ of a
free group $F_g \cong \pi_1(V^d_g)$ can be realized by the action of a group of
isometries of a hyperbolic handlebody $V^d_g$ (in the sense of the Nielsen realization
problem), for a sufficiently large dimension $d$.

\bigskip

{\bf Theorem 5.} {\sl  Let $G$ be a finite group of isometries of a
hyperbolic handlebody $V^d_g$, of dimension $d \ge 3$ and genus $g > 1$, which acts
faithfully on the fundamental group. 

\medskip

i) The order of $G$ is bounded by a polynomial of degree $d/2$ in $g$ if $d$ is even, 
and of degree $(d+1)/2$ if $d$ is odd. 

\medskip

ii) The degree $d/2$ is best possible in even dimensions whereas
in odd dimensions the optimal degree is  either $(d-1)/2$ or $(d+1)/2$.}

\bigskip

So in odd dimensions the optimal degree remains open at present; note  
that, for $d = 3$, the bound $(d+1)/2 = 2$ is not best possible since it gives a
quadratic bound instead of the actual linear bound $12(g-1)$, so maybe for all odd
dimensions the optimal degree is $(d-1)/2$.

\bigskip

{\bf 1.4}  \hskip 1mm   Finite group actions on finite graphs.

\medskip

Let $G$ be a finite group of automorphisms of a finite graph $\tilde\Gamma$ of rank $g>1$
(defined as the rank of its free fundamental group), allowing closed and multiple
edges.  Note that, without changing the rank of a 
graph, we can delete all {\it free edges}, i.e. nonclosed edges with a
vertex of valence 1 (an isolated vertex).  By possibly subdividing edges, we can also assume
that $G$ acts {\it without inversions} (of edges), i.e. no element acts on an edge as a
reflection in its midpoint.  We say that a finite graph is  {\it hyperbolic} if it has
rank $g>1$ and no free edges. In the following, all finite group actions on graphs will be
faithful and without inversions.

\medskip

By [Z2], each finite subgroup $G$ of the
outer automorphism group ${\rm Out} \, F_g$ of a free group $F_g$ can be induced by an
action of $G$ on a finite graph of rank $g$ (this is again a version of the Nielsen
realization problem).  Conversely, if $G$  acts on a hyperbolic graph then $G$ induces an
injection into ${\rm Out} \, F_g$   ([Z3, Lemma 1]). By [WZ], for $g \ge 3$ the
largest possible order of a finite subgroup of ${\rm Out} \, F_g$ is $2^gg!$ ([WZ]); in
particular, there is no linear or polynomial bound in $g$ for the order of $G$. 
In strong analogy with Theorems 1 and 4, the following holds (proved in section 3):

\bigskip
\vfill \eject

{\bf Theorem 6.}  {\sl i)  The maximal order of a finite
group $G$ acting with trivial edge stabilizers and without inversions
on a finite hyperbolic graph of rank $g$ is equal to
$6(g-1)$ or $4(g-1)$, and both cases occur for infinitely
many values of $g$.

\medskip

ii) Let $G$ be a finite group acting without inversions on a finite hyperbolic
graph. If $c$ denotes the order of an edge stabilizer of the action of $G$ then  $|G| \le
6c(g-1)$. 

\medskip

iii)  Equality $|G| = 6c(g-1)$ is obtained only for $c=1$, 2, 4, 8 and 16.    
There  are  infinitely many values of $c$ such that the second largest possibility
$|G| = 4c(g-1)$ is obtained.}

\bigskip

{\bf 2. Finite group actions on closed handles and the proof of Theorem 4}

\medskip

We briefly recall some concepts from [Z1]. 
Let $G$ be a finite group acting on a  handle
$H_g$. By the equivariant sphere theorem, there is an equivariant decomposition of $H_g$
into 0-handles $S^3$ connected by 1-handles $S^2 \times [-1,1]$, and an associated finite
graph $\tilde \Gamma$ with an action of $g$. In the language of [Z1], this induces on the
quotient orbifold $H_g/G$ the structure of a {\it closed handle-orbifold}, i.e. $H_g/G$ 
decomposes into  0-handle orbifolds 
$S^3/G_v$ connected by 1-handle orbifolds $(S^2/G_e) \times [-1,1]$ (in the case of a free
action, $H_g/G$ is a 3-manifold and one may just use the classical decomposition into prime
manifolds). This defines a finite graph of finite groups  $\G$ associated to the
$G$-action, with underlying graph 
$\Gamma = \tilde \Gamma/G$; by subdividing edges, we assume here that $G$ acts
without inversions of edges on $\Gamma$. The vertices of $\Gamma$ correspond to the
0-handle orbifolds, the edges to the 1-handle orbifolds. The vertex groups $G_v$ of $\G$ are
the stabilizers in $G$ of the 0-handles
$S^3$  of $H_g$ and isomorphic to finite subgroups of the orthogonal group SO(4), the edge
groups $G_e$ are stabilizers of 1-handles of $H_g$ and isomorphic to finite subgroups of
SO(3).  The fundamental group $\pi_1\G$ of the graph of groups $\G$ is defined as the
iterated free product with amalgamation and HNN-extension of the vertex groups along the
edge groups (starting with a maximal tree), and is isomorphic to the orbifold fundamental
group of the quotient orbifold $H_g/G$. There is a surjection 
$\phi: \pi_1\G \to G$, and $H_g$ is the orbifold covering of
$H_g/G$ associated to the kernel of $\phi$ (isomorphic to the free group
$F_g$). We will also assume in the following that the graph
of groups $\G$ has non {\it trivial edges}, i.e. edges with two different vertices such that
the edge group coincides with one of the two vertex groups.

\medskip

We denote by
$$\chi\G = \sum {1 \over |G_v|} - \sum {1 \over |G_e|}$$
the  Euler characteristic of the graph of groups $\G$ (the sum is taken over all
vertex groups $G_v$ resp. edge groups $G_e$ of $\G$); then 
$$g-1 =  -\chi\G \; |G|$$
(see [Se], [ScW] and [Z4] for the general theory of graphs of groups, groups acting on
trees and groups acting on finite graphs).

\bigskip

Remarks. i) The approach to finite group actions on 3-dimensional handlebodies
(Theorem 1) is analogous, using the equivariant Dehn lemma/loop theorem instead of the
equivariant sphere theorem. The 0-handles are disks
$D^3$ connected by 1-handles $D^2 \times [-1,1]$, the vertex groups of the graph of
groups $\G$ are finite subgroups of SO(3) and the edge groups finite subgroups of SO(2)
(i.e., cyclic groups).  In the case of maximal order $12(g-1)$,  $\pi_1\G$ is one of the
following four products with amalgamation ([Z2],[McMZ1]):

$$\D_2 *_{\Z_2} \S_3 , \hskip 4mm  \D_3*_{\Z_3}\A_4, \hskip 4mm  \D_4*_{\Z_4}\S_4,
\hskip 4mm  \D_5*_{\Z_5}\A_5$$ 

where $\D_n$ denotes the dihedral group of order $2n$.

\medskip

ii)  For the case of finite group actions on bounded surfaces (Theorem 2), one decomposes
the action along properly embedded arcs; the 0-handles are disks $D^2$ connected by
1-handles  $D^1 \times [-1,1]$, the vertex groups of $\G$ are
finite subgroups of O(2) (cyclic or dihedral) and the edge groups subgroups of ${\rm O}(1)
\cong \Z_2$ (i.e., of order two generated by a reflection of $D^1$, or trivial). In the
case of maximal order $12(g-1)$, $\pi_1\G$ is the free product with amalgamation
$\D_2 *_{\Z_2} \D_3$ (the first of the four groups in part i).

\bigskip

{\it Proof of Theorem {\rm 4}.}

\medskip

i)  Suppose that $G$ acts freely on a closed handle
$H_g$, $g > 1$. Since there are no orientation-preserving free actions of a finite group
on $S^2$, the edge groups of the associated graph of groups $\G$ are all trivial. It is easy
to see then that the minimum positive value for $-\chi\G$ is realized exactly by the graph
of groups $\G$ with exactly one edge and vertex groups $\Z_2$ and $\Z_3$, with $\pi_1\G
\cong \Z_2 * \Z_3$ and  $-\chi\G = 1-1/2 -1/3 = 1/6$ (we will say that $\G$ is of type (2,3)
in the following), and hence $|G| = 6(g-1)$ is the largest possible order.

\medskip

Let $M$ be the 3-manifold which is the connected sum of two lens spaces with fundamental
groups $\Z_2$ and $\Z_3$, with  $\pi_1(M) \cong \Z_2 * \Z_3$.  Let $\phi$ be a surjection
from $\Z_2 * \Z_3$ to a cyclic or dihedral group $G$ of order 6.
The regular  covering of $M$ associated to the kernel $F_2$ of $\phi$ is a closed handle of
genus 2 on which $G$ acts as the group of covering transformations, and this realizes the
largest possible order $|G| = 6(g-1)$. By factorizing $\Z_2 * \Z_3$ by characteristic
subgroups of $F_2$ of arbitrary large finite indices, one obtains examples for the maximal
order $|G| = 6(g-1)$ for arbitrarily large values of $g$ (see also the remark after the
proof for explicit examples realizing the maximum order).

\medskip

Concerning the lower bound $2(g+1)$, we consider a graph of groups $\G$ with exactly one
edge, of type $(2,n)$,   with $\pi_1\G \cong \Z_2 * \Z_n$ and $-\chi\G = 1-1/2 -1/n =
(n-2)/2n$. Let $M$ be the connected sum of two lens spaces, with  $\pi_1(M) \cong \Z_2 *
\Z_n$, and let 
$\phi: \Z_2 * \Z_n \to \D_n$ be a surjection onto the dihedral group of order $2n$.
The covering of $M$ associated to the kernel of $\phi$ is a closed handle $H_g$ of genus
$g$ with a free $G$-action; also, $g-1 = (n-2)|G|/2n =  n-2$, hence $n=g+1$ and $|G| =
2(g+1)$. (Note that, if $n$ is even, there is
also a surjection of $\Z_2 * \Z_n$ onto the cyclic group $\Z_n$ which gives an order $|G| =
n = 2g < 2(g+1)$.)

\medskip

It remains to show that $|G| = 2(g+1)$ is the largest possible order for infinitely many
values of $g$.  It is easy to see that the graphs of groups $\G$ with trivial edge
groups and with a possible surjection of $\pi_1\G$  onto a group of order $|G| > 2(g+1)$
have exactly one edge, of type $(2,n)$, (3,3), (3,4) or (3,5).

\medskip

We exclude first the case of a graph of groups $\G$ of type $(2,n)$. A finite
quotient $G$ of $\pi_1\G$ with torsionfree kernel has order $xn$, for some
positive integer $x$, hence $|G| = xn = 2n(g-1)/(n-2)$ and  $x(n-2) = 2(g-1)$.
Suppose that $g-1$ is a prime number. As seen above, the cases $x=1$ and $x=2$ give orders
$2g$ and $2(g+1)$, so we can assume that $x>2$ and hence $n=3$ or $n=4$.

\medskip

Let $n=3$, so $G$ has order $6(g-1)$. Suppose in addition that $g>7$;
then 6 and $g-1$ are coprime and, by a result of Schur-Zassenhaus, $G$ is a semidirect
product of
$\Z_{g-1}$ and a group $\bar G$ of order 6. We can also assume that 3 does not divide $g-2$
(or, equivalently, that $g-1$ is one of the infinitely many primes congruent to 2 mod 3,
by a result of Dirichlet). Then an element of order 3 in $\bar G$ acts 
trivially on $\Z_{g-1}$ by conjugation, the element of order 2 acts trivially or
dihedrally, and this implies easily that $G$ cannot be generated by two elements of orders
2 and 3.

\medskip

Now let $n=4$; then $G$ has order $4(g-1)$ and is a semidirect
product of $\Z_{g-1}$ and $\Z_4$ (since $g-1$ is prime). Suppose that $g-1$ is one of the
infinitely many primes congruent to 11 mod 12, so 4 does not divide
$g-2$ (and, as before, 3 does not divide $g-2$). Then the element of order two in $\Z_4$
acts trivially on $\Z_{g-1}$, is the unique element of order two in $G$ and 
the square of every element of order 4, so clearly $G$ cannot be generated by two elements
of orders 2 and 4.

\medskip

It remains to exclude the types (3,3), (3,4) and (3,5). If $\G$ is of type (3,3)
then $|G| = 3(g-1)$. As before, since 3 does not divide $g-2$, there is a unique subgroup
of order 3 in $G$ and $G$ is not generated by two elements of order 3. In the cases (3,4)
and (3,5) one has $|G| = 12(g-1)/5$ and $|G| = 15(g-1)/7$; excluding in addition $g=8$, also
these two cases are not possible.

\medskip

We have shown that  $m(g) = 2(g+1)$ for infinitely many genera $g$,
and this concludes the proof of part i) of Theorem 4.

\bigskip

ii)  For an odd integer $g$, we consider the semidirect product 
$G = (\Z_{(g-1)/2} \rtimes \Z_4) \rtimes \Z_2$, of order $4(g-1)$. Denoting by $x$ a
generator of 
$\Z_{(g-1)/2}$, by $y$ a generator of $\Z_4$ and by $t$ a generator of $\Z_2$, the actions
of the semidirect product are given by $\; yxy^{-1} = x^{-1}, \; txt^{-1} = x^{-1}$ and 
$tyt^{-1} = xy$.  There is a surjection with torsionfree kernel $\phi: \Z_2 * \Z_4 \to G$
which maps a generator of $\Z_2$ to $t$ and a generator of $\Z_4$ to $y$. As
before, $\phi$ defines a free action of $G$ on a closed handle $H_g$ of genus $g$, so 
$m(g) \ge 4(g-1)$; this leaves the possibilities $m(g) = 4(g-1)$
and $m(g) = 6(g-1)$.

\medskip

Suppose that $g = 2p+1$, for a prime $p>12$. We show that there is no surjection 
$\phi$ of $\Z_2 * \Z_3$ onto a group $G$ of order $6(g-1) = 12p$, and hence  
$m(g) = 4(g-1)$.  By the Sylow theorems,
such a group $G$ has a normal subgroup $\Z_p$, and the factor group is the alternating group
$\A_4$ (since this is the only group of order 12 generated by two elements of orders 2 and
3). Again by the theorem of Schur-Zassenhaus, $G$ is a semidirect product $\Z_p
\rtimes \A_4$. If the action of $\A_4$ on $\Z_p$ is trivial then clearly such a surjection
$\phi$ does not exist. Suppose that the action of $\A_4$ on $\Z_p$ is nontrivial; since the
automorphism group of
$\Z_p$ is cyclic, the action of $\A_4$ factors through a nontrivial action of the 
factor group $\Z_3$ of $\A_4$, and the subgroup $\D_2$ of $\A_4$ acts trivially. By
the Sylow theorems, up to conjugation we can assume that a surjection $\phi$ maps the factor
$\Z_3$ of $\Z_2 * \Z_3$ to $\A_4$. Since any involution in $\A_4$ acts trivially on $\Z_p$, 
every element of order 2 in $G$ is in $\A_4$, and hence $\phi$ is not 
surjective.

\medskip

Finally, any surjection $\phi: \Z_2 * \Z_3 \to \A_4$ defines a free action of $\A_4$ on a
closed handle of genus 3.  Factorizing by characteristic subgroups of arbitrary large
indices of the kernel $F_3$ of $\phi$, one obtains $m(g) = 6(g-1)$ for infinitely many
odd values of $g$.

\bigskip

iii)  Suppose that $n$ is prime, and let $g = n-1$ and $|G| =
2n$. Then it follows easily as above that $m(g) = 2n(g-1)/(n-2)$, hence infinitely many 
values of $n$ occur. We will show that also infinitely many values of $n$ do not occur.

\medskip

Let $n$ be congruent to 2 mod 8, $n \ne 2$,  and suppose that there exists $g$ such that
$m(g) = 2n(g-1)/(n-2)$. Then 8 divides  $(n-2)m(g) = 2n(g-1)$ and also $4(g-1)$, so 
$g-1$ is even and $g$ is odd. By ii), $m(g) = 2n(g-1)/(n-2) \ge 4(g-1)$, and this
gives the contradiction $n \le 4$.

\medskip

This concludes the proof of Theorem 4.

\bigskip

Remark.  We give an  explicit construction realizing the maximum order $|G| = 6(g-1)$
for infinitely many values of $g$.  Let  $g = p+1$, for a prime
$p$ such that 6 divides  $p-1$. We shall define  a surjection  $\phi$ of $\Z_2 * \Z_3$ onto
a semidirect product $G = \Z_p \rtimes \Z_6$;  this defines an action of $G$ on a
closed handle of genus $g = p+1$ of maximal possible order $|G|= 6(g-1)$, hence $m(g) =
6(g-1)$. Writing $\Z_p$ additively and
$\Z_6$ multiplicatively, suppose that a generator
$t$ of $\Z_6$ acts by conjugation on $\Z_p$ by an automorphism of order 6, in particular
$t^3$ acts dihedrally on $\Z_p$;  let
$\alpha$ be the automorphism of order 3 of $\Z_p$ induced by $t^2$.  Fixing a
generator
$a$ of $\Z_p$, one has
$\alpha(a) = a+b$ for some $b \ne 0$ in
$\Z_p$; then $\alpha^3(a) = a$ implies  $b + \alpha(b) + \alpha^2(b) = 0$. 
Considering the factors of $\Z_2 * \Z_3$, let $\phi$ map a generator of $\Z_2$ to $bt^3 =
t^3(-b)$ and a generator of $\Z_3$ to $bt^2$. Since also $b$
generates $\Z_p$, clearly $\phi$ is a surjection.

\bigskip 
\vfill \eject

{\bf 3. Finite group actions on finite graphs and the proof of Theorem 6}

\medskip

Let $G$ be a finite group acting without inversions on a finite, hyperbolic graph $\tilde
\Gamma$ of rank $g$.  Considering the quotient graph
$\Gamma = \tilde \Gamma/G$, we associate to each vertex group and edge group of
$\Gamma$ the stabilizer in $G$ of a preimage in
$\tilde \Gamma$ (starting with a lift of a maximal tree in $\Gamma$ to $\tilde \Gamma$);
this defines a finite graph of finite groups
$\G$ and a surjection $\phi: \pi_1\G \to G$, injective on vertex groups, with kernel
$F_g$.  Conversely, by the theory of groups acting on trees and graphs of groups (see 
[Se], [ScW] or [Z4]),  such a surjection $\phi: \pi_1\G \to G$  defines an action of
$G$ on a finite graph $\tilde \Gamma$ of rank $g =  -\chi\G \; |G| + 1$; the
action of $G$ on $\tilde \Gamma$ is faithful if and only if every finite normal subgroup of
$\pi_1\G$ is trivial (since a finite normal subgroup must be contained in all edge groups;
see [MeZ, Lemma 1]).

\bigskip

{\it Proof of Theorem {\rm 6}.}  

\medskip

i)  Let $G$ be a finite group which acts with trivial edge stabilizers and without
inversions on a finite graph $\tilde \Gamma$, of rank $g>1$. Then the associated graph of
groups $\G$ has trivial edge groups, and clearly $-\chi\G = 1/6$ is the smallest
positive value which can be obtained for the Euler characteristic $\chi\G$ (realized by the
graph of groups with one edge and edge groups $\Z_2$ and $\Z_3$).  Hence $|G| \le 6(g-1)$
and, as in the proof of Theorem 4, the upper bound $6(g-1)$ is obtained for infinitely many
values of $g$.

\medskip

For an integer $m>1$, choose a surjection $\phi:\Z_2 * \D_2 \to G$ where $G$ is the
dihedral group $\D_{2m}$ or the group $\Z_2 \times \D_m$, of order $4m$. Then
$\phi$ defines an action of $G$ on a finite graph of rank $g = m+1$, hence
$|G| = 4m = 4(g-1)$. On the other hand, if  $g-1$ is a prime
such that 3 does not divide
$g-2$ then it follows as in the proof of Theorem 4 that there does not exist a surjection of
$\Z_2 * \Z_3$ onto any group of order $6(g-1)$, so $4(g-1)$ is the maximal possible order
for infinitely many $g$.

\bigskip

ii) As before, the action of $G$ on $\tilde \Gamma$ is associated to a surjection with
torison-free kernel $\phi: \pi_1\G \to G$, for a finite graph of finite groups $\G$. We can
assume that
$\G$ has no trivial edges, i.e. edges with two different vertices such that the edge group
coincides with one of the vertex groups (by contracting such an edge). Since the action of
$G$ on $\tilde \Gamma$ is faithful, every finite normal subgroup of
$\pi_1\G$  is trivial. Let $e$ be an edge of $\Gamma$ with an edge group of order $c$; let
$\chi = \chi\G$ denote the Euler-characteristic of $\G$ and $n$ the order of
$G$.

\medskip

Suppose first that $e$ is a closed edge (a loop).  If $e$ is
the only edge of $\G$ then
$$-\chi  \ge  {1 \over c} - {1 \over 2c} = {1 \over c}, \hskip 5mm
g-1 = -\chi n \ge {n \over 2c},  \hskip 5mm  n \le 2c(g-1)$$
(since every finite normal subgroup of $\pi_1\G$ is trivial, the edge group of $e$ cannot
coincide with the vertex group).

\medskip

If $e$ is closed and not the only edge then
$$-\chi  \ge  {1 \over c}, \hskip 5mm
g-1 = -\chi n \ge {n \over c},  \hskip 5mm  n \le c(g-1).$$

Suppose that $e$ is not closed. If $e$ is the only edge of $\G$ then both
vertices of $e$ are isolated and
$$-\chi  \ge {1 \over c} - {1 \over 2c} - {1 \over 3c} = {1 \over 6c}, \hskip 5mm   g-1 =
-\chi \; n  \; \ge  \; {n \over 6c}, \hskip 5mm  n \le 6c(g-1).$$

\medskip

If $e$ is not closed, not the only edge and has exactly one isolated vertex then
$$-\chi  \ge {1 \over c} - {1 \over 2c} = {1 \over 2c}, \hskip 5mm   g-1 =
-\chi \; n  \; \ge  \; {n \over 2c}, \hskip 5mm   n \le 2c(g-1).$$
Finally, if $e$ is not closed, not the only edge and has no isolated vertex then
$$-\chi  \ge {1 \over c}, \hskip 5mm   g-1 =
-\chi \; n  \; \ge  \; {n \over c}, \hskip 5mm  n \le c(g-1).$$

\medskip

Concluding, in all cases we have $|G| \le 6c(g-1)$, proving ii).

\bigskip

iii)  By [G] and [DM], there are only finitely many free
products with amalgamation of two finite groups, without nontrivial finite normal
subgroups, such that the amalgamated subgroup has indices 2 and 3 in the two factors (and
the same holds also for indices 3 and 3). These  {\it effective {\rm
(2,3)}-amalgams} are classified in [DM], there are exactly seven such amalgams
(described below), and the amalgamated subgroups have order 1, 2, 4, 8 or 16. It follows
then from the proof of ii) that equality $|G| = 6c(g-1)$ can be obtained
only for these values of $c$.  

\medskip

On the other hand, by [D] there are infinitely many effective
(2,4)-amalgams, and hence $|G| = 4c(g-1)$ is obtained for infinitely many values of $c$.

\medskip

This concludes the proof of Theorem 6.

\bigskip

Remarks.  i) The seven effective $(2,3)$-amalgams, with amalgamated subgroups of orders $c =
1$, 2, 4, 8 or 16, are the following:
$$\Z_2 * \Z_3 , \hskip 4mm  \Z_4*_{\Z_2}\D_3, \hskip 4mm  \D_2*_{\Z_2}\D_3,
\hskip 4mm  \D_4*_{\D_2}\D_6, 
\hskip 4mm  \D_8 *_{\D_4} \S_4, \hskip 4mm \tilde \D_8 *_{\D_4} \S_4,$$
$$\Bbb  K_{32} *_{(\D_4 \times \Z_2)} (\S_4 \times \Z_2),$$
where $\tilde \D_8$ denotes the quasidihedral group of order 16 and $K_{32}$ a group of
order 32.

\medskip
\vfill \eject

ii) Finally, we describe the two
families which realize the largest possible orders for all $g$.
As noted before, the largest order of a finite group $G$ of automorphisms of a
finite graph of rank  $g>2$ without free edges (or equivalently, of a finite subgroup $G$
of ${\rm Out} \, F_g$) is  $2^gg!$, and this is realized by the automorphism
group $\; (\Z_2)^g \rtimes \S_g \;$ of a finite graph with one vertex and $g$ closed edges 
(a bouquet of $g$ circles or a  multiple closed edge), subdividing edges to avoid
inversions. Considering the quotient graph/graph of groups,  this action is
associated to a surjection
 
$$\phi:  \;  ((\Z_2)^g \rtimes \S_g) \; *_{((\Z_2)^{g-1} \rtimes \S_{g-1})} \;
((\Z_2)^g \rtimes \S_{g-1})  \; \to  \; (\Z_2)^g \rtimes \S_g$$

For $g \ge 3$, this realizes the unique action
of maximal possible order $2^gg!$ (the unique finite  subgroup of ${\rm Out} \, F_g$ of
maximal order, up to conjugation).

\medskip

The second family of large orders is given by the automorphism groups $\; S_{g+1} \times
\Z_2 \;$ of the graphs with two vertices and
$g+1$ connecting edges (a multiple nonclosed edge, subdividing edges again to avoid
inversions), associated to  the surjections 

$$\phi: \;  \S_{g+1} *_{\S_g} (\S_g \times \Z_2) \; \to \; 
\S_{g+1} \times \Z_2$$ 

These  realize the largest possible order for $g=2$, and again for
$g=3$.

\bigskip \bigskip

\centerline {\bf References}

\bigskip

\item {[CZ]} M.D.E. Conder, B. Zimmermann,  {\it Maximal bordered surface groups versus
maximal handlebodiy groups.}  Contemporary Math.  629 (2014), 99-105

\smallskip

\item {[D]} D. Djokovic,  {\it A class of finite group amalgams.} Proc. Amer. Math. Soc.
80 (1980), 22-26

\smallskip

\item {[DM]} D. Djokovic, G. Miller,  {\it Regular groups of automorphisms of cubic
graphs.}  J. Comb. Theory Series B 29 (1980), 195-230

\smallskip

\item {[G]} D.M. Goldschmidt,  {\it Automorphisms of trivalent graphs.}  Ann. Math. 111
(1980), 377-406

\smallskip

\item {[McMZ1]} D. McCullough, A. Miller, B. Zimmermann,  {\it Group actions on
handlebodies,}  Proc. London Math. Soc.  59   (1989), 373-415

\smallskip

\item {[McMZ2]} D. McCullough, A. Miller, B. Zimmermann,  {\it Group actions on nonclosed
2-manifolds.}  J. Pure Appl. Algebra  64   (1990),   269-292

\smallskip

\item {[MeZ]}  M. Mecchia, B. Zimmermann,  {\it  On finite groups of isometries of
handlebodies in arbitrary dimensions qnd finite extensions of Schottky groups.} 
Fund. Math.  230 (2015),  237-249

\smallskip

\item {[MiZ]} A. Miller, B. Zimmermann,  {\it  Large groups of symmetries of
handlebodies.}  Proc. Amer. Math. Soc. 106  (1989),  829-838

\smallskip

\item {[R]} F. Raymond,  {\it Classification of actions of the circle on 3-manifolds.}
Trans. Amer. Math. Soc. 131  (1968), 51-78

\smallskip

\item {[ScW]} P. Scott, T. Wall,  {\it  Topological methods in group theory.}
Homological Group Theory, London Math. Soc. Lecture Notes 36  (1979), Cambridge University
Press

\smallskip

\item {[Se]} J.P. Serre,  {\it  Trees.}  Springer, New York, 1980

\smallskip

\item {[WZ]} S. Wang, B. Zimmermann,  {\it  The maximum order finite groups of
outer automorphisms of free groups.}  Math. Z. 216  (1994), 83-87

\smallskip

\item {[Z1]}   B. Zimmermann,  {\it  On finite groups acting on a connected sum of
3-manifolds $S^2 \times S^1$.} Fund. Math. 226 (2014), 131-142

\smallskip

\item {[Z2]} B. Zimmermann,  {\it \"Uber Hom\"oomorphismen n-dimensionaler Henkelk\"orper
und endliche Erweiterungen von Schottky-Gruppen.}  Comm. Math. Helv. 56   (1981), 474-486

\smallskip

\item {[Z3]} B. Zimmermann,  {\it Finite groups of outer automorphism groups of free
groups,}  Glasgow Math. J. 38  (1996),  275-282

\smallskip

\item {[Z4]} B. Zimmermann,  {\it  Generators and relations for discontinuous groups.}
Generators and relations in Groups and Geometries,  NATO Advanced Study Institute Series
vol. 333  (1991), 407-436 (eds. Barlotti, Ellers, Plaumann, Strambach),  Kluwer Academic
Publishers

\smallskip

\item {[Z5]} B. Zimmermann,  {\it \"Uber Abbildungsklassen von
Henkelk\"orpern.}  Arch. Math. 33  (1979),  379-382

\bye